\documentclass{article}
\usepackage{amssymb,amsmath,theorem,euscript}

\newcounter{sec}

\def\sm{\smallskip}

\newcounter{punct}[sec]

\def\punct{\refstepcounter{punct}{\arabic{sec}.\arabic{punct}.  }}

\def\COUNTERS{\addtocounter{sec}{1}
              \setcounter{punct}{0}
          \setcounter{equation}{0}
          \setcounter{theorem}{0}
                  }

\newtheorem{theorem}{Theorem}[sec]
\newtheorem{proposition}[theorem]{Proposition}
\newtheorem{lemma}[theorem]{Lemma}

\begin{document}

 \def\ov{\overline}
\def\wt{\widetilde}
 \newcommand{\rk}{\mathop {\mathrm {rk}}\nolimits}
\newcommand{\Aut}{\mathop {\mathrm {Aut}}\nolimits}
\newcommand{\Out}{\mathop {\mathrm {Out}}\nolimits}
\renewcommand{\Re}{\mathop {\mathrm {Re}}\nolimits}
\def\Br{\mathrm {Br}}

\def\SL{\mathrm {SL}}
\def\SU{\mathrm {SU}}
\def\GL{\mathrm {GL}}
\def\U{\mathrm U}
\def\OO{\mathrm O}
 \def\Sp{\mathrm {Sp}}
 \def\SO{\mathrm {SO}}
\def\SOS{\mathrm {SO}^*}
 \def\Diff{\mathrm{Diff}}
 \def\Vect{\mathfrak{Vect}}
\def\PGL{\mathrm {PGL}}
\def\PU{\mathrm {PU}}
\def\PSL{\mathrm {PSL}}
\def\Symp{\mathrm{Symp}}
\def\End{\mathrm{End}}
\def\Mor{\mathrm{Mor}}
\def\Aut{\mathrm{Aut}}
 \def\PB{\mathrm{PB}}
 \def\cA{\mathcal A}
\def\cB{\mathcal B}
\def\cC{\mathcal C}
\def\cD{\mathcal D}
\def\cE{\mathcal E}
\def\cF{\mathcal F}
\def\cG{\mathcal G}
\def\cH{\mathcal H}
\def\cJ{\mathcal J}
\def\cI{\mathcal I}
\def\cK{\mathcal K}
 \def\cL{\mathcal L}
\def\cM{\mathcal M}
\def\cN{\mathcal N}
 \def\cO{\mathcal O}
\def\cP{\mathcal P}
\def\cQ{\mathcal Q}
\def\cR{\mathcal R}
\def\cS{\mathcal S}
\def\cT{\mathcal T}
\def\cU{\mathcal U}
\def\cV{\mathcal V}
 \def\cW{\mathcal W}
\def\cX{\mathcal X}
 \def\cY{\mathcal Y}
 \def\cZ{\mathcal Z}
\def\0{{\ov 0}}
 \def\1{{\ov 1}}
 \def\frA{\mathfrak A}
 \def\frB{\mathfrak B}
\def\frC{\mathfrak C}
\def\frD{\mathfrak D}
\def\frE{\mathfrak E}
\def\frF{\mathfrak F}
\def\frG{\mathfrak G}
\def\frH{\mathfrak H}
\def\frI{\mathfrak I}
 \def\frJ{\mathfrak J}
 \def\frK{\mathfrak K}
 \def\frL{\mathfrak L}
\def\frM{\mathfrak M}
 \def\frN{\mathfrak N} \def\frO{\mathfrak O} \def\frP{\mathfrak P} \def\frQ{\mathfrak Q} \def\frR{\mathfrak R}
 \def\frS{\mathfrak S} \def\frT{\mathfrak T} \def\frU{\mathfrak U} \def\frV{\mathfrak V} \def\frW{\mathfrak W}
 \def\frX{\mathfrak X} \def\frY{\mathfrak Y} \def\frZ{\mathfrak Z} \def\fra{\mathfrak a} \def\frb{\mathfrak b}
 \def\frc{\mathfrak c} \def\frd{\mathfrak d} \def\fre{\mathfrak e} \def\frf{\mathfrak f} \def\frg{\mathfrak g}
 \def\frh{\mathfrak h} \def\fri{\mathfrak i} \def\frj{\mathfrak j} \def\frk{\mathfrak k} \def\frl{\mathfrak l}
 \def\frm{\mathfrak m} \def\frn{\mathfrak n} \def\fro{\mathfrak o} \def\frp{\mathfrak p} \def\frq{\mathfrak q}
 \def\frr{\mathfrak r} \def\frs{\mathfrak s} \def\frt{\mathfrak t} \def\fru{\mathfrak u} \def\frv{\mathfrak v}
 \def\frw{\mathfrak w} \def\frx{\mathfrak x} \def\fry{\mathfrak y} \def\frz{\mathfrak z} \def\frsp{\mathfrak{sp}}
 \def\bfa{\mathbf a} \def\bfb{\mathbf b} \def\bfc{\mathbf c} \def\bfd{\mathbf d} \def\bfe{\mathbf e} \def\bff{\mathbf f}
 \def\bfg{\mathbf g} \def\bfh{\mathbf h} \def\bfi{\mathbf i} \def\bfj{\mathbf j} \def\bfk{\mathbf k} \def\bfl{\mathbf l}
 \def\bfm{\mathbf m} \def\bfn{\mathbf n} \def\bfo{\mathbf o} \def\bfp{\mathbf p} \def\bfq{\mathbf q} \def\bfr{\mathbf r}
 \def\bfs{\mathbf s} \def\bft{\mathbf t} \def\bfu{\mathbf u} \def\bfv{\mathbf v} \def\bfw{\mathbf w} \def\bfx{\mathbf x}
 \def\bfy{\mathbf y} \def\bfz{\mathbf z} \def\bfA{\mathbf A} \def\bfB{\mathbf B} \def\bfC{\mathbf C} \def\bfD{\mathbf D}
 \def\bfE{\mathbf E} \def\bfF{\mathbf F} \def\bfG{\mathbf G} \def\bfH{\mathbf H} \def\bfI{\mathbf I} \def\bfJ{\mathbf J}
 \def\bfK{\mathbf K} \def\bfL{\mathbf L} \def\bfM{\mathbf M} \def\bfN{\mathbf N} \def\bfO{\mathbf O} \def\bfP{\mathbf P}
 \def\bfQ{\mathbf Q} \def\bfR{\mathbf R} \def\bfS{\mathbf S} \def\bfT{\mathbf T} \def\bfU{\mathbf U} \def\bfV{\mathbf V}
 \def\bfW{\mathbf W} \def\bfX{\mathbf X} \def\bfY{\mathbf Y} \def\bfZ{\mathbf Z} \def\bfw{\mathbf w}
 \def\R {{\mathbb R }} \def\C {{\mathbb C }} \def\Z{{\mathbb Z}} \def\H{{\mathbb H}} \def\K{{\mathbb K}}
 \def\N{{\mathbb N}} \def\Q{{\mathbb Q}} \def\A{{\mathbb A}} \def\T{\mathbb T} \def\P{\mathbb P} \def\G{\mathbb G}
 \def\bbA{\mathbb A} \def\bbB{\mathbb B} \def\bbD{\mathbb D} \def\bbE{\mathbb E} \def\bbF{\mathbb F} \def\bbG{\mathbb G}
 \def\bbI{\mathbb I} \def\bbJ{\mathbb J} \def\bbL{\mathbb L} \def\bbM{\mathbb M} \def\bbN{\mathbb N} \def\bbO{\mathbb O}
 \def\bbP{\mathbb P} \def\bbQ{\mathbb Q} \def\bbS{\mathbb S} \def\bbT{\mathbb T} \def\bbU{\mathbb U} \def\bbV{\mathbb V}
 \def\bbW{\mathbb W} \def\bbX{\mathbb X} \def\bbY{\mathbb Y} \def\kappa{\varkappa} \def\epsilon{\varepsilon}
 \def\phi{\varphi} \def\le{\leqslant} \def\ge{\geqslant}

\def\UU{\bbU}
\def\Mat{\mathrm{Mat}}
\def\tto{\rightrightarrows}

\def\Gr{\mathrm{Gr}}

\def\graph{\mathrm{graph}}

\def\O{\mathrm{O}}

\def\la{\langle}
\def\ra{\rangle}

\begin{center}
\Large\bf
Sphericity \\ and multiplication of double cosets
\\ for infinite-dimensional classical groups

\bigskip

\large\sc
Yury A. Neretin\footnote{Supported by the grant FWF, P22122.}
\end{center}

{\small We construct spherical subgroups 
in infinite-dimensional classical groups 
$G$
(usually they are not symmetric and their finite-dimensional
analogs are not spherical). 
We present a structure of a semigroup on  double cosets
$L\setminus G/L$ for various subgroups  $L$ in $G$,
moreover these semigroups act in spaces of 
$L$-fixed vectors in unitary representations of
 $G$. We also obtain semigroup envelops of groups  $G$
 generalizing constructions of operator colligations.
}

\section{Introduction}

\COUNTERS

{\bf\punct Definition of spherical subgroups.\label{ss:spherical-def}}
 Let $G$ be a group,  $K\subset G$ a subgroup.
An irreducible unitary representation
 $\rho$  of $G$ in Hilbert space   $H$
is called {\it $K$-spherical},
if there exists a unique (up to a proportionality)
$K$-fixed vector  $v\in H$. 
The {\it spherical function} of a spherical representation
is given by the formula
$$
\Phi(g)=\la \rho(g) v,v\ra
.$$

Let $G$ be a semisimple Lie group,  $K$ a compact subgroup. The pair $(G,K)$ 
is called {\it spherical},
if any irreducible unitary representation
of  $G$ contains at most one non-zero
 $K$-fixed vector. The subgroup  $K$ is said to be a {\it spherical}
 subgroup of  $G$.

According to the Gelfand theorem
\cite{Gel}, symmetric subgroups are spherical.
In particular, maximal compact subgroups of semisimple groups are
spherical.

Kr\"amer classified \cite{Kra}
all spherical pairs   $(G,K)$ with simple group  $G$, 
the classification was extended to semisimple groups
in  \cite{Mik} and \cite{Bri}. 
All such pairs can be produced from symmetric pairs  $(G,K)$
by a minor diminishing
of the subgroup  $K$ or small enlargement of   $G$.
On spherical functions for non-symmetric pairs, see the work
of Knop \cite{Kno}.  

\sm

{\bf\punct Infinite-dimensional classical groups.%
\label{ss:icg0}} We consider Hilbert spaces 
 $l_2$ over $\R$, $\C$, or quaternions $\H$. 
 We call a matrix  $g$, whose elements are contained in $\R$, $\C$, $\H$,
{\it finite} if the matrix  $g-1$ has only finite number
of non-zero elements.
 We define {\it infinite-dimensional classical groups}
 $\U(\infty)$, $\OO(\infty)$, $\GL(\infty,\R)$,
 $\Sp(2\infty,\R)$, 
$\U(\infty,\infty)$, \dots as groups of finite matrices 
satisfying the usual identities.
All such groups  $G(\infty)$
are inductive limits of corresponding finite-dimensional
groups
 $G(n)$,
$$
\Sp(2\infty,\R)=\lim_{\longrightarrow} \Sp(2n,\R),\qquad
\U(\infty,\infty)=\lim_{\longrightarrow} \U(n_1,n_2),\qquad\dots
$$
We suppose that 
 $G(n)$ are canonically embedded to  $G(\infty)$ 
 and act to initial basis vectors. 
 If we consider groups of type  $\O(\infty,\infty)$ etc.
or products of groups, then we regard $n$ in $G(n)$
as an integer vector $n=(n_1,\dots)$.


{\bf\punct Heavy groups.%
\label{ss:icg}}
The orthogonal group $\O(\infty)$, the unitary group $\U(\infty)$,
and the quaternionic unitary group 
$\Sp(\infty)$ play the same role as compact groups in finite-dimensional
theory. We call them  {\it heavy}. 
For a heavy group 
 $K$, we denote by   $K^{\alpha}$
 the stabilizer of initial   $\alpha$
 basis elements in the corresponding  $\ell_2$.
 We write elements of  $K^{\alpha}$ as block 
$(\alpha+\infty)\times(\alpha+\infty)$ matrices of the form
 $\begin{pmatrix} 1&0\\0&h \end{pmatrix}$, where $h\in K$. 

 Let $\rho$ be a unitary representation of a heavy group
 in a Hilbert space
 $H$. Denote by  $H^{[\alpha]}$ the space of
  $K^{\alpha}$-fixed vectors.

\begin{theorem}
{\rm \cite{Olsh-fin}}
\label{th:KO}
The following conditions are
equivalent:
 
 --- $\rho$ is continuous with respect to the weak topology on  $K$;
 
 --- $\rho$ is continuous with respect to the uniform topology on  $K$;

--- $\cup H^{[\alpha]}$ is dense in $H$;

--- $\rho$ admits an extension to a weakly continuous 
representation of the semigroup of operators with norm
$\le 1$ in 
  $\ell_2$%
\footnote{We call such operators {\it contractions}.}.

If $\rho$ is irreducible, this is equivalent to  

--- $H^{[\alpha]}\ne 0$ for some $\alpha$.
\end{theorem}

We call such representations  {\it continuous.} 
Description of continuous representations is contained in
 \cite{Olsh-fin}, all such representations can be obtained by tensor operations
 from the simplest representation of
  $K$ in $\ell_2$.

Finite products of heavy groups also will be called 
 {\it heavy}. If
$K=K_1\times\dots\times K_p$, 
then we assume  $K^\alpha:=K_1^{\alpha_1}\times\dots\times K_p^{\alpha_p}$, 
now $\alpha=(\alpha_1,\dots,\alpha_p)$ is a multi-index. 
{\it Irreducible representations of a heavy 
group  $K=K_1\times\dots\times K_p$ are tensor products of 
representations of  $K_j$.}

\sm

{\bf\punct Olshanski's list.%
\label{ss:olsh}} The following list consists
of pairs  $(G,K)$= (classical group,  symmetric heavy subgroup). 
Explicits form of embeddings  $K(\infty)\to G(\infty)$ 
are sufficiently obvious.
\begin{align*}
&\text{Pairs of noncompact type} &&\text{Pairs of compact type} 
\\
& 
\bigl(\GL(\infty,\R), \O(\infty)\bigr) & & 
\bigl(\U(\infty), \O(\infty)\bigr)
\\
& 
\bigl(\GL(\infty,\C), \U(\infty) \bigr) &&  
\bigl( \U(\infty)\times \U(\infty), \U(\infty) \bigr)
\\
& 
\bigl(\GL(\infty,\H), \Sp(\infty) \bigr) & & 
\bigl(\U(2\infty),\Sp(\infty) \bigr)
\\&
\bigl(\U(\infty,\infty), \U(\infty)\times\U(\infty) \bigr) & & 
\bigl(\U(2\infty), \U(\infty)\times\U(\infty)  \bigr)
\\
& 
\bigl(\Sp(2\infty,\R),\U(\infty) \bigr) & & 
\bigl( \Sp(\infty),\U(\infty)\bigr)
\\
& 
\bigl(\Sp(2\infty,\C),\Sp(\infty) \bigr) & &
\bigl(\Sp(\infty) \times \Sp(\infty) ,\Sp(\infty)  \bigr)
\\
& 
\bigl(\Sp(\infty,\infty), \Sp(\infty)\times \Sp(\infty)   \bigr) && 
\bigl(\Sp(2\infty), \Sp(\infty)\times \Sp(\infty) \bigr)
\\&
\bigl(\O(\infty,\infty), \O(\infty)\times\O(\infty) \bigr) & & 
\bigl(\O(2\infty), \O(\infty)\times\O(\infty)  \bigr)
\\
& 
\bigl(\O(\infty,\C),\O(\infty), \bigr) && \bigl(\O(\infty)\times \O(\infty),\O(\infty) \bigr)
\\&
\bigl(\SOS(2\infty),\U(\infty) \bigr) && \bigl(\O(2\infty),\U(\infty) \bigr) 
\end{align*}

Also there are 3 series of symmetric pairs of finite rank:
\begin{align*}
&
\bigl(\U(p,\infty),\U(p)\times\U(\infty)\bigr)&&
\bigl(\U(p+\infty),\U(p)\times\U(\infty)\bigr)
\\
&
\bigl(\O(p,\infty),\O(p)\times\O(\infty)\bigr)&&
\bigl(\O(p+\infty),\O(p)\times\O(\infty)\bigr)
\\
&
\bigl(\Sp(p,\infty),\Sp(p)\times\Sp(\infty)\bigr)&&
\bigl(\Sp(p+\infty),\Sp(p)\times\Sp(\infty)\bigr)
\end{align*}

In pairs $(G,K)$ of the left column, the group 
$G$ ranges in all irreducible {\it infinite-dimensional classical groups.}%
\footnote{Of course, the groups 
 $\U(\infty)$, $\O(\infty)$, $\Sp(\infty)$ also are classical, $p=0$.}.
In the first table in the first column,
the subgroup $K$ is {\it a maximal heavy subgroup} in
 $G$.
The maximal heavy subgroup in   $\U(p,\infty)$ is $\U(\infty)$.

\sm

{\bf Sphericity theorem.} {\it All the above listed pairs  $(G,K)$ 
are spherical.}

\sm

Spherical representations are classified, see
 \cite{VK}, \cite{Olsh-dop}, \cite{OlshGB}, \cite{Pick},
 spherical functions are given by simple explicit expressions.
 For instance, consider the pair $\bigl(\GL(\infty,\R),\O(\infty)\bigr)$.
 Spherical representations are determined by finite collections
 of real parameters,
$s_1$, \dots, $s_p$ and also $a\in \R$ , $\sigma=0,1$.
Spherical functions are given by
$$
\Phi_{s,a,\sigma}(g)=
|\det g|^{ia}\cdot \mathrm{sgn}\cdot (\det g)^{\sigma} 
\prod_k\prod_l \left( \frac{1+is_k}2 \lambda_l+  \frac{1-is_k}2 \lambda_l^{-1}
\right)^{-1/2}
,$$
where
 $\lambda_j$ denote the singular numbers 
 of  $g\in \GL(\infty,\R)$.

{\sc Remark.} For the pair
$\bigl( \U(\infty)\times \U(\infty), \U(\infty) \bigr)$
a substantial harmonic analysis now is constructed, 
see \cite{BO},
apparently it exists for other pairs
of the second column,
see  \cite{Ner-hua}.\hfill $\square$

The class of spherical representations is not closed with respect to
natural operations (tensor products, restriction to subgroups),
it is reasonable to extend it.

\sm

{\bf\punct Admissible representations.%
\label{ss:admissible}} Fix a pair   $(G,K)$
from the list. We say that a unitary representation of
 $G$ {\it admissible} (\cite{OlshGB}), if it is continuous
 on the heavy subgroup  $K$. 
We also use as a synonym the term  "{\it representations of the pair $(G,K)$}".
Certainly, it is possible to define admissibility in terms 
of appropriate topologies on $G$ (see a discussion in \cite{Ner-book}).

\begin{lemma}
Spherical representations are admissible.
\end{lemma}

{\sc Proof.} It is easy to show that the subspace
 $\cup_\alpha H^{[\alpha]}$ is invariant with respect to 
$G$. On the other hand it is non-zero
 (it contains the spherical vector).
 In virtue of irreducibility, it is dense.
 \hfill$\square$

 \sm
 
 We emphasis that groups of finite
 matrices $\GL(\infty,\R)$, $\U(\infty)$, $\O(\infty,\C)$
 etc. are not type%
 \footnote{For definition, see, e.g., \cite{Dix}.} 
 $I$ groups, and theory of their unitary representations
 in the usual sense is impossible.
 The class of admissible representations is
 observable and substantial works 
 on representations of classical infinite-dimensional groups
 usually can be included to this scheme.
 
\sm

{\sc Remark.} The list contains pairs 
$\bigl(\O(2\infty), \O(\infty)\times\O(\infty)  \bigr)$, 
$\bigl(\O(2\infty), \U(\infty)  \bigr)$. The group $G$ 
in both cases is the same, but conditions of continuity
with respect to $K$ are different.
Therefore we get two different classes 
of representations. \hfill $\square$
 
\sm 
 
{\sc Remark.}  The list contains the pair  $(\U(\infty)\times\U(\infty),\U(\infty))$.
Its representations are not tensor products 
 $\rho_1\otimes\rho_2$ of representations of 
$\U(\infty)$ (the corresponding finite-dimensional
theorem, see  \cite{Dix}, 13.11.8,
in this case is not valid), 
on representations of this pair,
see
 \cite{SV}, \cite{VK}, \cite{OlshGB}. \hfill $\square$

\sm


{\bf\punct Product of double cosets and train.%
\label{ss:train}}
Let  $\rho$  be an admissible representation
of the pair $(G,K)$ in the space  $H$.
As above, let
 $H^{[\alpha]}$ be spaces of  $K^\alpha$-fixed
 vectors. 
 Denote by  $P^{[\alpha]}$ the operator of orthogonal projection
 to  $H^{[\alpha]}$. We define operators 
 $$
 \ov\rho_{\alpha,\beta}(g) = P^{[\beta]} \rho(g): H^{[\alpha]}\to H^{[\beta]}
 .$$
 It easy to verify that for any
 $h_1\in K^\alpha$, $h_2\in K^\beta$ we have
$$
\ov\rho_{\alpha,\beta}(h_2 g  h_1)=\ov\rho_{\alpha,\beta}(g).
$$
Therefore the function $\ov\rho_{\alpha,\beta}$ 
is defined on double cosets
$\frg\in K^\beta\setminus G/K^\alpha$.

\sm

{\bf Multiplicativity theorem.} {\it There exists
a natural multiplication
$(\frg, \frh)\mapsto \frg\circ \frh$,
$$
K^\gamma\setminus G/ K^\beta \,\, \times\,\, K^\beta\setminus G/K^\alpha
\,\,\to\,\,
 K^\gamma\setminus G/K^\alpha
,$$
defined for all
 $\alpha$, $\beta$, $\gamma\in\Z_+$. Moreover, for
 any admissible unitary representation
 $\rho$ of the pair $(G,K)$ the following equality holds
\begin{equation}
 \ov\rho_{\beta,\gamma}(\frg) \ov\rho_{\alpha,\beta}(\frh)= \ov\rho_{\alpha,\gamma}(\frg\circ\frh)
 \label{eq:repcat}
.\end{equation}
}

We get a category, we call it by
 {\it train} $\frT(G,K)$ of the pair $(G,K)$. 
Objects of the train are indices (multi-indices)
 $\alpha\ge 0$, and morphisms are double cosets.
  The formula (\ref{eq:repcat})
claims that we get a representation of a category 
(see \cite{Ner-book}, II.5),
i.e., a functor from  $\frT(G,K)$
to the category of Hilbert spaces and linear bounded operators.

The operation
$g\mapsto g^{-1}$ induces  the map
$K^\alpha\setminus G/ K^\beta \to
K^\beta\setminus G/ K^\alpha$, we denote it by
 $\frg\mapsto \frg^*$.
It satisfies the identity
$$
(\frg\circ\frh)^*=\frh^*\circ \frg^*
.$$
The representation  $\ov\rho$ is a {\it $*$-representation}
in the following sense:
$$
\rho(\frg^*)=\rho(\frg).
$$

{\bf Approximation property.} {\it Let   $\tau$
be a representation of the category
$\frT(G,K)$.
Let  $\|\tau(\frg)\|\le 1$ for all $\frg$, $\tau(1)=1$.
Then the representation  $\tau$ has the form $\ov\rho$
for uniquely defined representation  $\rho$  of the pair $(G,K)$.}


\sm

{\bf\punct Comparison with finite-dimensional
case.%
\label{ss:Hecke}} Let $G$ be a Lie group,  $L$
a compact subgroup, not necessarily maximal.
Denote by  $\cM(L\setminus G/L)$ the space of finite measures that are invariant
with respect to left and right shifts by
elements of $L$. Evidently, $\cM(L\setminus G/L)$ 
is an algebra with respect to the convolution.
For instance, if  $G$ is a  $p$-adic  $\GL(n)$ and $L$ 
iwahoric subgroup, we get the Hecke algebra
 (see. \cite{Mac}); if $G$ is a real Lie group of rank 1
and  $L$ the maximal compact subgroup, then
we get the hypergroup of generalized translate, see \cite{Koo}.

 Let $\rho$ be a unitary representation of  $G$ 
 in the space $H$. Denote by 
$H^L$ the space of  $L$-fixed vectors. For $\mu\in \cM(L\setminus G/L)$
the operator $\rho(\mu)$ has the block form 
$\begin{pmatrix}A(\mu)&0\\0&0 \end{pmatrix}$ 
with respect to the decomposition $H=H^L\oplus (H^L)^\bot$.
Thus we get a functor from unitary representations of
 $G$ to representations of the convolution algebra 
$\cM(L\setminus G/L)$.

It turns out that for infinite-dimensional groups 
the convolution of double cosets can degenerate to 
a multiplication. A first example. apparently was discovered by R.S.Ismagilov \cite{Ism1},
\cite{Ism2},
further \cite{Olsh-tree}, \cite{Olsh-fin}, \cite{OlshGB}, \cite{Ner-sto},
numerous examples are discussed in \cite{Ner-book}.

\sm


{\bf \punct Explicit formula for the product.%
\label{ss:circ}} See \cite{OlshGB}.
For definiteness, consider the pair 
$(G,K)=(\GL(\infty,\R),\O(\infty))$ and double cosets
$K^\beta\setminus G/ K^\alpha$. 
In other words, we consider matrices 
 $\frg=\begin{pmatrix}A&B\\C&D\end{pmatrix}$ of  size $(\beta+\infty)\times(\alpha+\infty)$
 up to the equivalence 
\begin{equation}
  \begin{pmatrix}A&B\\C&D\end{pmatrix} \sim
 \begin{pmatrix}1_\beta&0\\0&U\end{pmatrix} 
\begin{pmatrix}A&B\\C&D\end{pmatrix} \begin{pmatrix}1_\alpha&0\\0&V\end{pmatrix}
\label{eq:colligation-1}
, \end{equation}
where $U$, $V$ are orthogonal matrices, $1_\alpha$ 
denotes the unit matrix of size 
$\alpha$. The product in the category  $\frT$
is given by  
\begin{equation} 
\begin{pmatrix}A&B\\C&D\end{pmatrix}
 \circ \begin{pmatrix}P&Q\\R&T\end{pmatrix}
 := \begin{pmatrix}A&B&0\\C&D&0\\0&0&1\end{pmatrix}
 \begin{pmatrix}P&0&Q\\0&1&0\\R&0&T\end{pmatrix} =
\begin{pmatrix} AP&B&AQ\\ CP&D&CQ\\ R&0&T \end{pmatrix}
\label{eq:product}
\end{equation}

Note that a similar (but not precisely the same) algebraic structure is known
as a product of operator colligations, see, e.g.,  \cite{Bro}.

\sm

{\bf\punct Characteristic functions%
\label{ss:characteristic}} (below we do not develop this topic).
It appears  (see \cite{Ner-book}, IX.4), that  $\circ$-multiplication
can be transformed to a more customary operation%
\footnote{This is an analog of Livshits characteristic function, which is
a tool in spectral theory of non-self-adjoint operators,
see \cite{Liv}, \cite{Pot}, \cite{Bro}. The Livshits function also  take part
in the representation theory of the pair 
 $\bigl(\U(1+\infty)\times\U(\infty),\U(\infty)\bigr)$,
see. \cite{Olsh-CR}.}
in the following way. Fix $\lambda\in\C\cup \infty$. 
Write the equation
$$
\begin{pmatrix}
p_+\\ \lambda x\\ p_-\\x
\end{pmatrix}
=
\begin{pmatrix} 
  \begin{matrix}A&B\\C&D\end{matrix}&0
  \\
  0&  \begin{pmatrix}A&B\\C&D\end{pmatrix}^{t-1}
\end{pmatrix}
\begin{pmatrix}
q_+\\ y\\q_-\\ \lambda y
\end{pmatrix}
.$$
Consider all vectors  $(p_+,p_-)\oplus(q_+, q_-)\in\C^{2\beta}\oplus\C^{2\alpha}$,
for which this equation has a solution as an equation with indeterminates
$x$, $y$. We get a subspace  
$\chi(\lambda)$ of half-dimension in
$\C^{2\beta}\oplus\C^{2\alpha}$ depending on $\lambda$.
We regard it as an relation (a multi-valued map)
from
 $\C^{2\alpha}$ to $\C^{2\beta}$.

On this language, $\circ$-multiplication
corresponds to the product
of characteristic functions, i.e., pointwise
product of relations.

Note that
 $\chi(\lambda)$
 is a rational map from Riemann sphere to the Lagrangian  Grassmannian
 in
 $\C^{2\alpha}\oplus\C^{2\beta}$.
The function
$\chi(\lambda)$ is also  $J$-contractive in the Potapov sense
\cite{Pot}.


\sm

{\bf\punct Self-similarity.%
\label{ss:self-similarity}} See. \cite{Olsh-symm}. 
Let $(G,K)=\bigl(\GL(\infty,\R), \O(\infty)\bigr)$.
Denote by $G^\alpha\subset G$ 
the stabilizer of initial  $\alpha$ basis vectors in
$\ell_2$. Then the pair $(G^\alpha,K^\alpha)$ is isomorphic to  $(G,K)$.
Similar pairs are defined in all the cases, but a definition requires some care:
we take basis vectors fixed by
 $K^\alpha$,  then $G^\alpha$ is their stabilizer in $G$.
 
 \sm

{\it Fix an unitary representation  $\rho$
of  $(G,K)$. 
Its restriction to a sufficiently small subgroup
 $(G^\alpha,K^\alpha)\simeq (G,K)$ contains a spherical subrepresentation.
 All spherical representations of $(G,K)$ obtained in this way
 are equivalent {\rm(}for all $\alpha${\rm)}.}

\sm

{\bf\punct Mantle.%
\label{ss:mantle}} Let $(G,K)=\bigl(\GL(\infty,\R), \O(\infty)\bigr)$.
Split $\N$ into two countable subsets  $\N=\Xi\cup\Omega$
(for instance, odd and even numbers).
Denote by  $G^\Xi\subset G$ 
(respectively $G(\Xi)$) the stabilizer of basis vectors
enumerated by elements of $\Xi$ (respectively, $\Omega$).
In other words, $G^\Xi$ and $G(\Xi)$ consist 
of matrices of the forms 
$
\begin{pmatrix}
1&0\\0& g
\end{pmatrix}$
and
$
\begin{pmatrix}
g&0\\0& 1
\end{pmatrix}
$, 
where blocks correspond to the decomposition  $\ell_2(\N)=\ell_2(\Xi)\oplus\ell_2(\Omega)$.
Evidently, $\bigl(G(\Xi), K(\Xi)\bigr)$ is isomorphic
to  $(G,K)$. We fix the isomorphism 
  
 $$
 i_\tau:(G,K)\to\bigl(G(\Xi), K(\Xi)\bigr)
 ,$$
 induced by some bijection $\tau:\N\to \Xi$.

 Let $\rho$ be a unitary representation of 
$(G,K)$ in the space $H$. Denote
by  $H^{[\Xi]}$ the space  of
$K^\Xi$-fixed vectors. By $P^{[\Xi]}$ we denote the operator of projection
on
 $H^{[\Xi]}$.

\sm

{\bf Lemma about self-restriction.}
{\it The space $H^{[\Xi]}$ is 
$G(\Xi)$-invariant. The representation  $\rho\circ i_\tau$
of the pair  $(G,K)$
 in $H^{[\Xi]}$ is equivalent to 
$\rho$.}

\sm

On the other hand, for each  $g\in G$ we assign the operator
 $\ov\rho_{\Xi,\Xi}(g):H^{[\Xi]}\to H^{[\Xi]}$
as 
$$\ov\rho_{\Xi,\Xi}(g)=P^{[\Xi]} \rho(g).$$
As above, $\ov\rho$
is a function on double cosets.

\sm

{\bf Version of multiplicativity theorem.} {\it 
Double cosets 
$\Gamma_\infty=K^\Xi\setminus G/K^\Xi$ admit a natural structure
of a semigroup 
{\rm(}{\bf mantle} of the pair $(G,K)${\rm)},
the map $\ov\rho_{\Xi,\Xi}$ is a representation of
 $\Gamma_\infty$.}
 
 \sm

Product is given by the same formula 
(\ref{eq:product}), in the present situation 
all blocks of the matrix are infinite.

Next, notice that the natural map
$G(\Xi)\simeq G \to K^\Xi\setminus G/K^\Xi$
is injective, therefore we observe that the representation
of the semigroup 
$K^\Xi\setminus G/K^\Xi$ extends the representation $\rho$ of $G(\Xi)\simeq G$.

As a result, we get

\sm

{\bf Theorem about mantle.} {\it Any unitary
representation of the pair 
$(G,K)$ can be extended canonically
to a representation of the semigroup 
 $\Gamma_\infty$.}

\sm

{\bf \punct Spherical characters.\label{ss:character}}
Notice, that  $Z:=K^\Xi\setminus G^\Xi/K^\Xi\subset K^\Xi\setminus G/K^\Xi$
is a central subsemigroup in the mantle.
In any irreducible representation 
$\rho$ of the pair
 $(G,K)$ it acts by scalar operators.
 This determines a multiplicative character
 $\Phi_\rho:Z\to \C^\times$.

If $\rho$ is spherical, then $\Phi_\rho$
coincides with the spherical function of $\rho$.
Generally, it coincides with spherical functions of subgroups as in
 \ref{ss:self-similarity}.

\sm


{\bf \punct Non-standard pairs.%
\label{ss:Ness}} Until this moment we discussed pairs $(G,K)$
containing in the Olshanski list. Now let
 $G=\GL(\infty,\R)\times\dots\times\GL(\infty,\R)$,
let $K\subset G$ be the group $\O(\infty)$ embedded to $G$ by the diagonal.

\sm

a) The pair $(G,K)$ is spherical. Notice, that this statement has no
finite-dimensional analogs. Multiplicativity theorem also holds.
 
\sm

b) Nessonov \cite{Ness1}, \cite{Ness2} obtained description
of all spherical representations
of this pair.

\sm

c) The paper \cite{Ner-char} contains a construction
of characteristic functions, it  is valid in our case.

\sm

There arises the following question:
In which generality the statements of Subsections 
\ref{ss:olsh}--\ref{ss:character} hold?

\sm

{\bf \punct Purposes of the paper.}
We introduce class of {\it $(G,L)$-pairs},
where $G$ is a classical group, $L$ 
is a heavy subgroup. For this class the following  hold true:

\sm

--- the multiplicativity theorem; 

\sm

--- the approximation property;

\sm

---  the construction of the mantle; 

\sm

--- existence of spherical characters.

\sm

The construction of characteristic functions
given in \cite{Ner-char} survives in these cases.
  (in the present paper this topic is not discussed).
  
  \sm

Also, we define a subclass of  {\it pure $(G,L)$-pairs}.
In this case there hold

--- sphericity theorem;

--- self-similarity.

Definition of  $(G,L)$-pairs is contained in Section 2,
sketches of proofs in Section 3,  Section 4 contains conjectures
and additional remarks. We do not repeate formulations of theorems.


 \section{Definition}

\COUNTERS

{\bf\punct Straight embeddings of heavy groups..%
\label{ss:right}} First, let heavy groups $L$, $K$
be irreducible, i.e, they have the form 
 $\U(\infty)$, $\O(\infty)$, $\Sp(\infty)$. 
 There are obvious homomorphisms:
\begin{align*}
&\O(\infty)\to \U(\infty) && \O(\infty)\to\Sp(\infty)&& 
\U(\infty)\to\Sp(\infty)
\\
&
\U(\infty)\to\O(2\infty)&&\Sp(\infty)\to \O(4\infty)&&
\Sp(\infty)\to\U(2\infty)
.\end{align*}
Also element-wise complex conjugation determines
a homomorphism
 $\U(\infty)$ to itself.
We call such homomorphisms and also identical homomorphisms
 $g\mapsto g$
{\it trivial.}

Next, let  $K$ be irreducible, 
$L=L_1\times \dots\times L_p$. We say that a homomorphism 
$L\to K$ is {\it straight}, if it has the form
\begin{equation}
(g_1,\dots,g_p)\mapsto
\begin{pmatrix}
1_a& 0&\dots &0
\\
0& \tau_1(g_{m_1})&\dots&0
\\
\vdots&\vdots &\ddots &\vdots
\\
0&0&\dots&\tau_N(g_{m_N})
\end{pmatrix},\qquad
\text{$g_k\in L_k$}
\label{eq:emb}
,\end{equation}
where $\tau_j$ are trivial homomorphisms.
We admit unit matrix  $1_a$ of  arbitrary order 0, 1, 2, \dots, $\infty$. 
The number
$N$
can be in the limits
 $0<N<\infty$. We do not require a homomorphism to be faithful.

Finally, let  $L=L_1\times \dots\times L_p$,  $K=K_1\times \dots\times K_q$.
Then a {\it straight homomorphism} is a product of  straight homomorphisms
$L\to K_1$, \dots, $L\to K_q$.

\sm


{\bf\punct $(G,L)$-pairs.%
\label{ss:GL}} We say that a {\it classical group}
$G=G_1\times \dots\times G_r$ is a finite product 
of irreducible classical groups
(see the list above).  
Let $K=K_1 \times \dots\times K_q$ be the maximal
heavy subgroup
in  $G$.
Let $L$ be a straightly embedded heavy
subgroup in  $K$.
We call such objects  {\it $(G,L)$-pairs.}

A pair $(G,L)$ is {\it pure}, if the centralizer
of the subgroup  $L$
in $G$ is trivial%
\footnote{Examples. The pair
$(G,L)=\bigl(\GL(3\infty,\R), \O(\infty)\times \O(\infty)\times \O(\infty)\bigr)$ is pure.
The centralizer of  $L$ in the group of all bounded
operators is $(\R^\times)^3$.
But it is not contained in the group  $G=\GL(3\infty,\R)$
of finite matrices. The pair
$\bigl(\O(1+3\infty,\R), \O(\infty)\times \O(\infty)\times \O(\infty)\bigr)$
is not pure, the centralizer in $\O(1+3\infty,\R)$ is  $\Z_2$.}.
An equivalent condition: in matrices  (\ref{eq:emb})
the unit block  $1_a$ is absent and there are no groups
$\U(p,\infty)$, $\O(p,\infty)$, $\Sp(p,\infty)$
among factors $G_j$%
\footnote{This limitation is necessary, because
for such groups the maximal heavy subgroup is not a precise analog
of the maximal compact group.}.

\sm

We say that a pair is {\it finite}, 
if  blocks $1_a$ in formula (\ref{eq:emb})
for all factors $K_i$ have finite sizes.
Finite pairs will appear only in conjectures in the last section.

\sm


{\bf\punct Definition of multiplication on double cosets.%
\label{ss:def-product}}
Let $L$ be irreducible. Fix  $\alpha\ge 0$ and consider
the sequence 
$\Theta_m^{[\alpha]}\in L$ given by
\begin{equation}
\rho(\Theta_m^{[\alpha]})=
\begin{pmatrix}
1_\alpha&0&0&0
\\
0&0&1_m&0
\\
0&1_m&0&0
\\
0&0&0&1_\infty
\end{pmatrix}
\label{eq:Theta}
.
\end{equation}
If $L=L_1\times\dots\times L_p$ is reducible, we fix an multi-index     $\alpha$ and take the net depending on  a multi-index
$m$,
$$
\Theta_{m}^{[\alpha]}= \bigl(\Theta_{m_1}^{[\alpha_1]},\dots
\Theta_{m_p}^{[\alpha_p]}
\bigr).
$$

Let
$$
\frg\in K^{[\gamma]}\setminus G/K^{[\beta]}
,\qquad
\frh\in K^{[\beta]}\setminus G/K^{[\alpha]}
.
$$
Choose representatives
 $g\in\frg$, $h\in\frh$.

Consider subsequence 
$$z_m=g\Theta_m^{[\beta]}h\in G
.
$$

\begin{proposition}
a{\rm)} A double coset  $\frz_m=K^{[\gamma]}z_m K^{[\alpha]}$
is eventually constant.

\sm

b{\rm)} The result $g\circ h=\lim \frz_m$ 
does not depend on a choice of representatives 
 $g\in\frg$, $h\in\frh$.

\sm

c{\rm)} The operation obtained in this way
is associative, i.e., for each
 $\alpha$, $\beta$, $\gamma$, $\delta$ and any
$$
\frf\in
K^{\delta}\setminus G/K^{\gamma}
,\qquad
\frg\in K^{\gamma}\setminus G/K^{\beta}
,\qquad
\frh\in K^{\beta}\setminus G/K^{\alpha}
$$
we have $(\frf\circ\frg)\circ\frh=\frf\circ(\frg\circ\frh)$.
\end{proposition}

As a result we get a category,
we call it {\it train} of the pair $(G,L)$ and denote
by
$\frT(G,L)$.

Instead of formal proof in general case, we present examples explaining why this happens.


\sm

{\bf\punct Example.%
\label{ss:example1}} Consider the pair $\bigl(\GL(\infty,\R),\O(\infty)\bigr)$.
To simplify notation, assume $\gamma=\beta=\alpha$. Fix
$g$, $h\in G$. Write  $g$, $h$ as block matrices of the size
 $\alpha+N+\infty$.
 Let $N$ be sufficiently large. Then   $g$, $h$ have the form
 \begin{equation}
 g=\begin{pmatrix}a&b&0\\ c&d&0\\0&0&1_\infty \end{pmatrix}\qquad
  h=\begin{pmatrix}p&q&0\\ r&t&0\\0&0&1_\infty \end{pmatrix}
 .
 \label{eq:gh}
 \end{equation}
The product $g\Theta^{[\alpha]}_{N+k}h$ equals
$$
\begin{pmatrix}
ap& aq & 0 & b & 0 & 0
\\
cp & cq & 0 & d &0 & 0
\\
0 & 0 & 0 & 0& 1_k&0
\\
r & t &0 &0 &0 &0
\\
0& 0 & 1_k & 0 & 0 &0
\\
0&0& 0&0 &0&  1_\infty
\end{pmatrix}
.$$
Since we examine double cosets, we transpose rows and columns 
whose numbers 
$>\alpha$. In this way we can reduce our matrix to arbitrary
of two convenient forms:
$$
S_1=
\begin{pmatrix}
ap& aq & b& 0
\\
cp& cq& d&0
\\
r&t& 0& 0
\\
0&0&0&1_{2k+\infty}
\end{pmatrix}
\qquad
S_2=
\begin{pmatrix}
ap&b& aq &  0
\\
cp&d& cq& 0
\\
r&0&t&  0
\\
0&0&0&1_{2k+\infty}
\end{pmatrix}
.$$
We observe that the result is independent on
 $k$.
 Next, we must verify independence on a choice of  representatives, i.e.,
 $gu$, $vh$ c $u$, $v\in K_\alpha$ 
 produce the same double coset.
 Without loss of generality, we can assume 
$u$, $v\in\O(N)$ (otherwise we choose larger $N$ earlier). 
We replace
$g$, $h$ from (\ref{eq:gh}) to
 \begin{equation}
 g=\begin{pmatrix}a&bu&0\\ c&du&0\\0&0&1 \end{pmatrix}\qquad
  h=\begin{pmatrix}p&q&0\\ vr&vt&0\\0&0&1\end{pmatrix}
 .
 \end{equation}
with orthogonal matrices  $u$, $v$.
The  new matrix   $S_2$ is
$$
\begin{pmatrix}
ap&bu& aq &  0
\\
cp&du& cq& 0
\\
vr&0&vt&  0
\\
0&0&0&1
\end{pmatrix}
=
\begin{pmatrix}
1&&&
\\
& 1&&
\\
&&v&
\\
&&&1
\end{pmatrix}
\begin{pmatrix}
ap&b& aq &  0
\\
cp&d& cq& 0
\\
r&0&t&  0
\\
0&0&0&1
\end{pmatrix}
\begin{pmatrix}
1&&&
\\
& u&&
\\
&&1&
\\
&&&1
\end{pmatrix}
.$$

Associativity is clear from the  formula
 (\ref{eq:product}).

\sm

{\bf\punct Two examples.%
\label{ss:example2}} a) Consider the pair
$G=\GL(\infty, \R)\times\dots\times\GL(\infty,\R)$
with the  diagonal subgroup  $L=\O(\infty)$. 
We have a finite collection of matrices 
$g_j=
\begin{pmatrix}
a_j&b_j\\c_j&d_j
\end{pmatrix}
$ 
defined up to the equivalence
\begin{multline*}
\left\{
\begin{pmatrix}
a_1&b_1\\c_1&d_1
\end{pmatrix},
\dots,
\begin{pmatrix}
a_q&b_q\\c_q&d_q
\end{pmatrix}
\right\}
\sim\\ \sim
\left\{
\begin{pmatrix}
1_\beta&0
\\
0&u
\end{pmatrix}
\begin{pmatrix}
a_1&b_1\\c_1&d_1
\end{pmatrix}
\begin{pmatrix}
1_\alpha&0
\\
0&v
\end{pmatrix}
,\dots,
\begin{pmatrix}
1_\beta&0
\\
0&u
\end{pmatrix}
\begin{pmatrix}
a_q&b_q\\c_q&d_q
\end{pmatrix}
\begin{pmatrix}
1_\alpha&0
\\
0&v
\end{pmatrix}
\right\}
.\end{multline*}

Then  $\circ$-product is defined component-wise by 
the formula
(\ref{eq:product}).
 
\sm

b) Let $G=\U(2\infty)$, $K=\O(\infty)\times\O(\infty)$.
Now we have unitary matrices  $g$
determined up to the equivalence 
{\small
$$
\begin{pmatrix}
a_{11}&b_{11}& a_{12}& b_{12}
\\
c_{11}&d_{11}& c_{12}& d_{12}
\\
a_{21}&b_{21}& a_{22}& b_{22}
\\
c_{21}&d_{21}& c_{22}& d_{22}
\end{pmatrix}
\sim
\begin{pmatrix}
1&&&
\\
& u_1&&
\\
&&1&
\\
&&&u_2
\end{pmatrix}
\begin{pmatrix}
a_{11}&b_{11}& a_{12}& b_{12}
\\
c_{11}&d_{11}& c_{12}& d_{12}
\\
a_{21}&b_{21}& a_{22}& b_{22}
\\
c_{21}&d_{21}& c_{22}& d_{22}
\end{pmatrix}
\begin{pmatrix}
1&&&
\\
& v_1&&
\\
&&1&
\\
&&&v_2
\end{pmatrix}
,$$
}
where $u_1$, $u_2$, $v_1$, $v_2\in\O(\infty)$.

The product is given by
\begin{equation}
g\circ g'=
\begin{pmatrix}
a_{11}&b_{11}&0& a_{12}& b_{12}&0
\\
c_{11}&d_{11}&0& c_{12}& d_{12}&0
\\
0& 0& 1& 0& 0& 0
\\
a_{21}&b_{21}&0& a_{22}& b_{22}&0
\\
c_{21}&d_{21}&0& c_{22}& d_{22}&0
\\
0& 0& 0& 0& 0& 1
\end{pmatrix}
\begin{pmatrix}
a'_{11}&0& b'_{11}& a'_{12}&0& b'_{12}
\\
0& 1& 0& 0& 0& 0
\\
c'_{11}&0& d'_{11}& c'_{12}&0& d'_{12}
\\
a'_{21}&0& b'_{21}& a'_{22}&0& b'_{22}
\\
0& 0& 0& 0& 1& 0
\\
c'_{21}&0&d'_{21}& c'_{22}&0& d'_{22}
\label{eq:shtrih}
\end{pmatrix}
.
\end{equation}

Double cosets are taken with respect to the group of matrices of
the form
$$
\begin{pmatrix}
1 &&&&&
\\
&u_{11}&u_{12}&&&
\\
&u_{21}&u_{22}&&&
\\
&&&1&&
\\
&&&& v_{11}&v_{12}
\\
&&&& v_{21}&v_{22}
\end{pmatrix},
\qquad
\text{where\,}
\begin{pmatrix}
u_{11}&u_{12}
\\
u_{21}&u_{22}
\end{pmatrix},
\begin{pmatrix}
v_{11}&v_{12}
\\
v_{21}&v_{22}
\end{pmatrix}
\in \O(2\infty)
.
$$

\section{Proofs}

\COUNTERS

{\bf\punct Multiplicativity%
\label{proof-multiplicativity}}
(the analog of  \ref{ss:train}).

\begin{lemma}
\label{l:basic}
For any unitary representation
  $\pi$ of a heavy group  $L$
the sequence $\pi(\Theta_m^{[\alpha]})$ {\rm(}see {\rm(\ref{eq:Theta}))}
weakly converges to the projector  $P^{[\alpha]}$ to $H^{[\alpha]}$.
\end{lemma}

See \cite{Ner-book}, Theorem VIII.1.4. 
A priory proof is given
for a symmetric group  $S(\infty)$,
 proofs for  $\O(\infty)$, $\U(\infty)$, $\Sp(\infty)$
are the same. On the other hand, it is easy to
reduce the statement from the explicit classification 
of representations of heavy groups \cite{Olsh-fin}, see also \cite{Ner-book}. \hfill $\square$

\sm

Let us prove the multiplicativity theorem for an arbitrary pair
 $(G,L)$.
Let $\frg$, $\frh$, $g$, $h$  be the same as in Subsection
 2.3. Then
$$
\ov \rho(\frg \circ \frh)=P^{[\gamma]}\rho(g\Theta_m^{[\alpha]}h): H^{[\alpha]}
\to H^{[\gamma]}
$$
for sufficiently large
 $m$. On the other hand, 
\begin{multline*}
\lim_{m\to\infty}
P^{[\gamma]}\rho(g\Theta_m^{[\alpha]}h)=
\lim_{m\to\infty}
P^{[\gamma]}\rho(g)\rho(\Theta_m^{[\alpha]})\rho(h)
=\\=
P^{[\gamma]}\rho(g)\bigl(\lim_{m\to\infty}\rho(\Theta_m^{[\alpha]})\bigr)\rho(h)
=P^{[\gamma]}\rho(g) P^{[\beta]}\rho(h)
= \ov\rho(\frg)\ov\rho(\frh)
.
\qquad \square
\end{multline*}

Thus, for any unitary representation of
 $(G,L)$ we construct a representation
 of the category  $\frT(G,L)$. 
 Now we present the inverse construction.


\sm

{\bf\punct Approximation property.%
\label{ss:proof-approximation}}
Let $\alpha<\beta$ (therefore, $L^\alpha\supset L^\beta$). 
Denote by 
 $\lambda_{\alpha,\beta}$ the double coset  
$L^\beta\cdot 1\cdot L^\alpha$. Set $\mu_{\beta,\alpha}:= \lambda_{\alpha,\beta}^*$.
Then the following identities hold:
\begin{equation}
\mu_{\beta,\alpha}\circ\lambda_{\alpha,\beta}=1^{[\alpha]},
\qquad
\lambda_{\beta,\gamma}\lambda_{\alpha,\beta}=\lambda_{\alpha,\gamma}
,
\label{eq:ordered}
\end{equation}
where $1^{[\alpha]}$ denotes the unit  $L^\alpha\cdot 1\cdot L^\alpha$
of the semigroup $L^\alpha\setminus G/ L^\alpha$. 
We get a structure of  
{\it ordered category} in the sense of  \cite{Ner-book}, III.4.

Notice that $\psi_{\alpha,\beta}:=\lambda_{\alpha,\beta}\circ\mu_{\beta,\alpha}$ as a subset in
 $G$ coincides with $L^\alpha$. But it is a nontrivial idempotent
 in 
$L^\beta\setminus G/ L^\beta$:
$$
\psi_{\alpha,\beta}^2=\psi_{\alpha,\beta}\qquad\psi^*=\psi_{\alpha,\beta}
.$$

Consider a  $*$-representation  $R$ of the train  $\frT(G,L)$
by contractive operators,
let $R$ assigns unit operators to units of semigroups of endomorphisms.
Denote by  $H^{[\alpha]}$ the Hilbert space,
corresponding to a (multi)index $\alpha$.
In virtue of (\ref{eq:ordered}), an operator $R(\lambda_{\alpha,\beta})$ is an operator of isometric embedding  
 $H^{[\alpha]}\to H^{[\beta]}$. The  projection operator
 to the image of the embedding is
$\psi_{\alpha,\beta}$.
Now we have a chain of isometric
embedding  (to be definite, let  $\alpha=(\alpha_1,\alpha_2)$
be a bi-index):
$$
\begin{matrix}
\dots&        &                       &        &H^{\alpha_1+1,\alpha_2} &   && &\dots \\
     &\searrow&                       &\nearrow&                & \searrow &   &\nearrow \\                  
\dots&        & H^{\alpha_1,\alpha_2} &        & &          &H^{\alpha_1+1,\alpha_2+1}& &\dots\\
     &\nearrow&                       &\searrow&                & \nearrow &   &\searrow \\                  
\dots&        &                       &        &H^{\alpha_1,\alpha_2+1} &  && &\dots
\end{matrix}
$$

Denote
by $H$ the limit space.

The group $G$ is a union of finite-dimensional
groups
 $G(\alpha)$.
 Notice that the natural map 
 $G(\alpha)$ to
$L^\alpha\setminus G/ L^\alpha$ is an isomorphism.
Therefore we get representation of
$G(\alpha)$ in the space $H^{[\alpha]}$. 
Operators of representation are contractive with
their inverses. Therefore they are unitary.

Hence the group $G$ acts in $H$. We omit a watching of the remaining details.
\hfill $\square$


\sm

{\bf\punct Sphericity.%
\label{proof:spherical}} 

\begin{lemma}
For a pure $(G,L)$-pairs the semigroup  $L\setminus G/L$
is commutative.
\end{lemma}

The case is that product in such
$L\setminus G/L$ means that we "join" one matrix to another. For instance, formula 
 (\ref{eq:shtrih}) in this case is reduced to 
$$
g\circ g'=
\begin{pmatrix}
d_{11}&d_{12}
\\
d_{21}&d_{22}
\end{pmatrix}
\circ
\begin{pmatrix}
d'_{11}&d'_{12}
\\
d'_{21}&d'_{22}
\end{pmatrix}
=
\begin{pmatrix}
d_{11}&0&d_{12}&0
\\
0&d'_{11}&0&d'_{12}
\\
d_{21}&0&d_{22}&0
\\
0&d'_{21}&0&d'_{22}
\end{pmatrix}
.$$
To get $g'\circ g$, we conjugate right-hand side 
by a matrix
$$
\begin{pmatrix}
J&0\\0&J
\end{pmatrix},\qquad
\text{where}\qquad
J=
\begin{pmatrix}
0&0&1_N&0\\
0&1_\infty&0\\
1_N&0&0&0\\
0&0&0&1_\infty
\end{pmatrix}
$$
with sufficiently large  $N$. \hfill $\square$


Thus $L\setminus G/L$  is an Abelian semigroup with
involution. Its irreducible representations are one-dimensional.
On the other hand, irreducibility
of representation of  $(G,L)$
is equivalent to irreducibility of the corresponding
representation of the train
 $\frT(G,L)$;
 but for an irreducible representation of an ordered category,
 all representations of semigroups of endomorphisms are irreducible
 (see \cite{Ner-book}, Lemma III.4.3).
 This implies also the following statement.

\begin{proposition}
Spherical functions of the pair $(G,L)$ 
are homomorphisms%
\footnote{non-arbitrary.}
from the semigroup $L\setminus G/L$ to the multiplicative
group  $\C^\times$.
\end{proposition}

\sm


{\bf\punct Self-similarity.%
\label{proof:self-similarity}} 
Consider a pure pair $(G,L)$ and its subgroup 
$G^\alpha$ as in 1.8.

\begin{lemma}
\label{l:center}
The image of $G^\alpha$ in $L^\alpha \setminus G/L^\alpha$ 
is a central semigroup canonically isomorphic 
to $L\setminus G/L$.
 \end{lemma}

{\sc Proof.} The isomorphism is induced by
 the "shift"  
$G\to G^\alpha$, i.e., $g\mapsto \begin{pmatrix}1&0\\0&g \end{pmatrix}$.
Further, let  $g\in G$, $h\in G^\alpha$. 
It is easy to see that for sufficiently large
 $m$,
$$
g\,(\Theta_m^{[\alpha]} h \Theta_m^{[\alpha]})
=
(\Theta_m^{[\alpha]} h \Theta_m^{[\alpha]})\,g
\quad
\Rightarrow
\quad
L^{\alpha}\cdot g \Theta_m^{[\alpha]}h\cdot L^{\alpha}
=
L^{\alpha}\cdot h \Theta_m^{[\alpha]}g \cdot  L^{\alpha},
\qquad\square
$$

Next, for
 $\alpha<\beta$
there is  a natural map $K^\beta\setminus G/K^\beta\to
K^\alpha\setminus G/K^\alpha$ that can be defined by the formula
$$
\Pi:\frg\mapsto  \mu_{\beta,\alpha} \circ \frg\circ \lambda_{\alpha,\beta}
.$$

\begin{lemma}
The  projection $\Pi:K^\beta\setminus G^\beta/K^\beta\to
K^\alpha\setminus G^\alpha/K^\alpha$ 
of central semigroups coincides with isomorphism induced
by the shift 
 $G^\beta\to G^\alpha$.
\end{lemma}

{\sc Proof.} The element  $g\in G^\alpha$ 
and the corresponding element
$\wt g\in G^\beta$ are conjugate in  $G^\alpha$ 
by an element of the group $K^\alpha$. \hfill $\square$

\sm

Therefore the semigroup
 $L\setminus G/L$ acts in the whole space 
$H$. If a representation  $\rho$ is irreducible,
then this semigroup acts by scalar operators.
The corresponding character coincides  with spherical functions
of groups  $(G^\alpha,K^\alpha)$ acting in $H$. 

\sm


{\bf \punct  Self-restriction.%
\label{proof:self-restriction}} 
Let us prove the statement for
$(G,L)=\bigl(\GL(\infty,\R),\O(\infty))$.
We preserve the notation of \ref{ss:mantle}.

The map
 $\tau:\N\to\Xi$ determines  0-1 matrix,
 denote it by $A$,
 the norm of $A$ is 1.
 According Theorem \ref{th:KO}, representations
 of a heavy group admits a weakly continuous  extension
 to representation of  semigroup of contructions.
In particular, for any continuous representation
$\pi$ of $L$ in a Hilbert space $H$, there is an operator $\pi(A)$ assigned
to the matrix  $A$. 
 
\begin{lemma}
 $\pi(A)$ is an isometric map from $H$ 
 to the subspace 
 $H^{[\Xi]}$.
\end{lemma}

{\sc Proof.} 
First, $(\pi(A))^*\pi(A)=\pi(A^*A)=\pi(1)=1$. Therefore 
$\pi(A)$ is an isometric embedding. Similarly,  
$\pi(A)(\pi(A))^*=\pi(AA^*)$.  The operator  $AA^*$ 
is the projection to the subspace 
$\ell_2(\Xi)$ in $\ell_2(\N)$.
Next, we again refer to explicit form of representations of heavy groups,
 \cite{Olsh-fin}.\hfill$\square$

\sm

Consider an admissible representation $\rho$ of a pair $(G,L)$.
 Denote by  $\pi$ the restriction 
 of 
$\rho$ to $L$.
 Extend $\pi$ to the semigroup of contractions.  
 
 \begin{lemma}
 The operator $\pi(A)$ intertwines the representation 
$\pi$  of the pair $(G,L)$ in the space $H$
and the representation
$\pi\circ i_\tau$ in the space   $H^{[\Xi]}$.
 \end{lemma}
 
 {\sc Proof.}
Consider a sequence of matrices-permutations   $B_{m}$ weakly  converging to
$A$.
Then $\pi(B_m)$  weakly converges to $\pi(A)$.
The group  $G$ is an inductive limit of 
reductive groups  $G(n)$. For each group $G(n)$
the operator $\pi(B_m)$ with sufficiently large number
 $m$ is intertwining. Therefore the limit is intertwining
 for the whole group $G$.\hfill $\square$ 

\sm

 Let $L=L_1\times \dots\times L_p$. Each factor is a group
 of unitary operators in a certain space
 $\ell_1=\ell_2(\N)$.
 We split each copy of
 $\N$ into two countable subsets $\N=\Xi_j\cup\Omega_j$.
Fix bijections $\tau_j:\N\to\Xi_j$. 
Now we can repeat the same words.


 \sm
 
 {\bf\punct Mantle.} Theorem on mantle follows from multiplicativity theorem
 and lemma about self-restriction.

\sm

 
 {\bf \punct Spherical characters,} see \ref{ss:character}. 
 Consider a $(G,L)$-pair, acting in a direct sum
of the spaces  $\ell_2$ as in \ref{ss:right}--\ref{ss:GL}.
Consider all basis vectors
 $e_\mu$ in $\oplus \ell_2$ that are fixed by the group
 $L$.
Consider the subgroup $G_{min}\subset G$ stabilizing all $e_\mu$.
 
 \sm
 
 {\sc Example.} Let $(G,L)= \bigl(\GL(7+2\infty,\R),\O(\infty)\times\O(\infty)\bigr)$.
 Then $G_{min}=\GL(2\infty,\R)$.
 
 \sm
 
 \begin{proposition}
 Mantle of the pair $(G,L)$ contains a central semigroup
 canonically isomorphic 
$L\setminus G_{min}/L$.
 \end{proposition}
 
 Explain this by an example. Let $(G,L)=\bigl(\U(p+\infty)\times\U(\infty),\U(\infty)\bigr)$,
 the subgroup $L\simeq\U(\infty)$ is embedded to the product
 by the diagonal,
 $G_{min}=\U(\infty)\times\U(\infty)$.
 Redenote $(G,L)$ as $\bigl(\U(p+2\infty)\times\U(2\infty),\U(2\infty)\bigr) $,
 we write its elements as pairs of block matrices
  $(g,h)$,
 where the size of $g$ is $(p+\infty+\infty)$ and the size of
 $h$ is $(\infty+\infty)$.
 Consider a subgroup $Z\simeq\U(\infty)\times\U(\infty)$
 consisting of pairs of matrices
 $$
 \left\{
 \begin{pmatrix}1&&\\&1&\\&&d_1 \end{pmatrix}, \begin{pmatrix}1&\\&d_2 \end{pmatrix}\right\}
 ,$$
 where $d_1$, $d_2\in\U(\infty)$.
 Consider a subgroup     $M\subset Z$
 consisting of pairs satisfying $d_1=d_2$.
 Then the mantle is the semigroup  $M\setminus G/M$. 
 It contains the central subgroup  $M\setminus Z/M$.
 Property of centrality is a version of
 Lemma
\ref{l:center} \hfill $\square$


 \section{Discussion}
 
\COUNTERS 
 
 {\bf\punct General conjectures.}
 Let a pair $(G,L)$ be finite, see  .\ref{ss:GL}.
Then

\sm

a) Unitary representations of $(G,L)$ have type $I$.

\sm

b) There exists only countable set of representations
with a given spherical character.

\sm

c) Let $(G,L)\to (G',L')$ be an embedding 
of pairs such that the embedding 
$L\to L'$ is straight in the sense of \ref{ss:right}.
Then the restriction of any irreducible
representation $\rho$ of the pair  $(G',L')$ to
 $(G,L)$ has a pure discrete spectrum.  In particular, this is the case
 for tensor products
 of irreducible representations.

\sm
 
d) An explicit classification 
of spherical representations is possible.
Precisely, they can be obtained by the construction
described in  \cite{Ner-book}, IX.5 (the difference between 
symmetric and non-symmetric pairs is explained in
 IX.5.6).

e) In \cite{Ner-char} for any element $\frg\in L^\alpha\setminus G/L^\beta$ 
there were obtained spectral data: "characteristic function"
(it is an inner function of a matrix argument)
and "eigensurface" (a hypersurface in a Grassmannian).
The element
$\frg$ can be uniquely reconstructed from these data.


\sm

{\bf\punct Why we consider only such pairs?}
There are many embeddings of $\U(n)$ 
to unitary groups  $\U(N)$ of big size.
For instance, $g\mapsto g\otimes g$ determines an embedding 
$\U(n)\to \U(n^2)$.
But for finite  $g\in\U(\infty)$ the matrix  $g\otimes g$
is not finite. 
I.e., formally,  embeddings $\U(n)$ to
 $\U(\cdot)$, different from straight embeddings do not survive
 as  $n\to\infty$.

Certainly, this argument is not sufficient. We can search 
natural subgroups in the complete operator group in
 $\ell_2(\N\times\N)$ containing matrices  
$g\otimes g$. The author do not see such subgroups.%
\footnote{For instance, the subgroup generated by finite unitary matrices
and matrices  $g\otimes g$ is a direct product of two generating subgroups.}. 

\sm

{\bf\punct  Sphericity.} There are limits of spherical
pairs that are not pure, say
$\bigl(\U(1+\infty)\times \U(\infty),\U(\infty)\bigr)$,
$\bigl(\O(1+2\infty), \U(\infty)\bigr)$, see \cite{Olsh-CR}.

Also there are many cases of sphericity with respect to subgroups of the form
 $L\times M$,
where  $L$ is  heavy  and  $M$  is compact
(see, e.g., the pair  $\U(p+\infty)/\U(p)\times\U(\infty)$ mentioned above).

\sm


{\bf \punct Symmetric groups.}
For symmetric groups similar facts hold, an example is discussed
in  \cite{Ner-symm}. Our proofs are based on Lemma \ref{l:basic}, 
they survive literally for symmetric groups.

{\tt Math.Dept., University of Vienna,

 Nordbergstrasse, 15,
Vienna, Austria

\&

Institute for Theoretical and Experimental Physics,

Bolshaya Cheremushkinskaya, 25, Moscow 117259,
Russia

\&

Mech.Math.Dept., Moscow State University,

Vorob'evy Gory, Moscow

e-mail: neretin(at) mccme.ru

URL:www.mat.univie.ac.at/$\sim$neretin

wwwth.itep.ru/$\sim$neretin
}

\begin{thebibliography}{cc}

\bibitem{BO}
Borodin, A., Olshanski, G. {\it Harmonic analysis on the infinite-dimensional unitary group and determinantal point processes.}  Ann. of Math. (2)  161  (2005),  no. 3,
1319--1422. 

\bibitem{Bri}
 Brion, M. {\it Classification des espaces homog\`enes sph\'eriques}, Compositio Math. 63 (1987), 189--208.

\bibitem{Bro}
Brodski, M. S.,
{\it Unitary operator colligations and their characteristic functions.}
Russ. Math. Surveys, 1978, 33:4, 159--191.

\bibitem{Gel}
Gelfand, I. M.
{\it Spherical functions in symmetric Riemann spaces.} (Russian)
Doklady Akad. Nauk SSSR (N.S.) 70, (1950). 5--8. English transl. in
Transl., Ser. 2, Amer. Math. Soc. 37 (1964), 39--43.


\bibitem{Dix}  Dixmier, J.
{\it Les $C^*$-alg\`ebres et leurs repr\'esentations.}
Gauthier Villars, 1964. 


\bibitem{Ism1}
Ismagilov, R.S.,
{\it Elementary spherical functions on the groups $\SL(2,P)$ over
a field $P$, which is not locally compsct with respect to
the subgroup of matrices with integral elements.} 
Math. USSR-Izvestiya, 1967, 1:2, 349--380. 


\bibitem{Ism2}
Ismagilov, R.S., {\it Spherical functions over a normed field whose residue 
field is infinite.} 
  Funct. Anal. and Its Appl.
Volume 4, N. 1, 37--45.




\bibitem{Kno}
Knop, F. {\it Semisymmetric polynomials and the invariant theory of matrix vector pairs.}
  Represent. Theory  5  (2001), 224--266.

\bibitem{Koo}
 Koornwinder, T. H. 
{\it Jacobi functions and analysis on noncompact semisimple Lie groups.}
In
{\it  Special functions: group theoretical aspects and applications,}
  1--85, Math. Appl., Reidel, Dordrecht, 1984.

\bibitem{Kra}
M. Kr\"amer, {\it Sph\"arische Untergruppen in kompakten zusammenhe\"angenden Liegruppen},
 Compositio Math. 38 (1979), no. 2, 129--153.


\bibitem{Liv}
Livshits, M.S. {\it On a certain class of linear operators in Hilbert space.}
 Mat. Sb., N. Ser. 19(61), 239-262 (1946),
English transl. in
Amer. Math. Soc. Transl. (Ser. 2), Vol. 13, 61--83
 (1960)


\bibitem{Mac}
Macdonald, I. G. {\it Symmetric functions and Hall polynomials.}
 Second edition. With contributions by A. Zelevinsky. 
  Oxford University Press, New York, 1995.

\bibitem{Mik}
Mikityuk, I. V.  {\it On the integrability of invariant Hamiltonian systems 
with homogeneous
 configuration spaces},
Math. USSR-Sb. 57 (1987), 527--546.

\bibitem{Ner-sto}
Neretin, Yu. A.
{\it 
Categories of bistochastic measures and representations of some infinite-dimensional groups.}
Russian Acad. Sci. Sb. Math. 75 (1993), no. 1, 197--219. 

\bibitem{Ner-book}
   Neretin, Yu. A. {\it Categories of symmetries and
 infinite-dimensional groups.}
 Oxford University Press, New York, 1996.

\bibitem{Ner-hua}
 Neretin, Yu. A. {\it Hua type integrals over unitary groups 
 and over projective limits of unitary groups,} Duke Math. J. 114 (2002), 239--266.

\bibitem{Ner-symm}
Neretin, Yu. A., {\it Infinite tri-symmetric group, multiplication 
of double cosets, and checker topological field theories.}
Preprint, available via {\tt arXiv:0909.4739}.

\bibitem{Ner-char} 
Neretin, Yu. A. {\it Multi-operator colligations and multivariate spherical functions.}
Preprint, available via {\tt arXive:1006.2275}.

\bibitem{Ness1}
 Nessonov, N. I.,   {\it Factor-representation of the group $GL(\infty)$
and admissible representations $\GL(\infty)^X$}, (Russian)  Math. Phys., Analysis,
Geometry (Kharkov Math. Journ.), 2003, N4,  167--187.

\bibitem{Ness2}
Nessonov, N.I.
{\it Factor-representation of the group $\text{$GL$}(\infty)$
and admissible representations of $\text{$GL$}(\infty)^X$. II.}
(Russian)
Mat. Fiz. Anal. Geom. (Kharkov Math. Journ.) 2003, No. 4, 524--556.


\bibitem{Olsh-tree}
Olshanski, G. I.
{\it New ``large'' groups of type ${\rm I}$.}
  Current problems in mathematics, Vol. 16 (Russian),  pp. 31--52, 228, Akad. Nauk SSSR,
 Vsesoyuz. Inst. Nauchn. i Tekhn. Informatsii (VINITI), Moscow, 1980.
English transl. J. Sov. Math. 18 (1982) 22--39.

\bibitem{Olsh-short}
Olshanski, G. I.
{\it  Unitary representations of infinite-dimensional pairs $(G,\,K)$
 and the formalism of R. Howe.}
 Soviet Math. Dokl. 27 (1982), no. 2, 290--294.

\bibitem{Olsh-fin}
Olshanski, G.I
{\it Infinite-dimensional groups of finite $R$-rank: description of representations 
and asymptotic theory.}
 Funct. Anal. and Appl., 1984, 18:1, 22--34.

\bibitem{Olsh-dop}
Olshanski, G. I.
{\it Unitary representations of the group ${\rm SO}_0(\infty,\infty)$
 as limits of unitary representations of the groups ${\rm SO}_0(n,\infty)$ as $n\to \infty$}. Funct. Anal. Appl. 20 (1986), no. 4, 292--301. 


\bibitem{OlshGB}
Olshanski, G.I.
{\it Unitary representations of infinite dimensional pairs
$(G,K)$ and the formalism of R. Howe.}
In
{\it Representation of Lie groups and related topics,}
Adv. Stud. Contemp. Math. 7, 269--463 (1990).

\bibitem{Olsh-symm}
 Olshanski, G.I.,
{\it Unitary representations of $(G,K)$-pairs connected with
the infinite symmetric group $S(\infty)$.}
Leningr. Math. J. 1, No.4, 983--1014 (1990).

\bibitem{Olsh-CR} Olshanski, G.I.
{\it Caract\`eres g\'en\'eralis\'es du groupe $U(\infty)$ et fonctions int\`erieures.} Comptes Rendus
Acad. Sci. Paris. S\'er. 1, 313 (1991), 9-12.

\bibitem{Pick}
Pickrell, D. {\it  Separable representations for automorphism 
groups of infinite symmetric spaces.}  J. Funct. Anal.  90  (1990),  no. 1, 1--26.
 
\bibitem{Pot}
Potapov, V. P. {\it The multiplicative structure of $J$-contractive
 matrix functions.} Trudy Moskov. Mat. Obschestva. 4 (1955), 125--236;
English transl. 
Amer. Math. Soc. Transl. (2)  15  (1960), 131--243.

\bibitem{SV}
Str\u atil\u a, ┼., Voiculescu, D.,
 {\it Representations of AF-algebras and of the group $U(\infty )$}.
 Lecture Notes in Mathematics, Vol. 486 (1975).

\bibitem{VK}
Vershik, A. M.; Kerov, S. V. Characters and factor-representations 
of the infinite unitary  group. 
 Soviet Math. Dokl. 26, 570--574 (1983).



\end{thebibliography}
\end{document}